\documentclass[12pt]{amsart}
\usepackage{amssymb}
\usepackage[dvips]{graphics}
\ifx\mathrm\undefined\let\mathrm\rm\fi
\ifx\mathbf\undefined\let\mathbf\bf\fi
\ifx\mathfrak\undefined\let\mathfrak\frak\fi
\ifx\mathcal\undefined\let\mathcal\cal\fi
\ifx\mathbb\undefined\let\mathbb\Bbb\fi
\ifx\emph\undefined\let\emph\it\fi
\newcommand{\g}{{{\mathfrak g}\,}}

\newcommand{\Z}{{\mathbb Z}}

\newcommand{\R}{{\mathbb R}}

\newcommand{\Ref}[1]{{(\ref{#1})}}
\newcommand{\be}{\begin{displaymath}}
\newcommand{\ee}{\end{displaymath}}
\newcommand{\bea}{\begin{eqnarray*}}
\newcommand{\eea}{\end{eqnarray*}}

\newcommand{\dontprint}[1]{\relax}

\newtheorem%
{thm}{Theorem}[section]
\newtheorem%
{proposition}[thm]{Proposition}
\newtheorem%
{lemma}[thm]{Lemma}
\newtheorem%
{lemmadef}[thm]{Lemma-Definition}
\newtheorem%
{corollary}[thm]{Corollary}
\newtheorem%
{conjecture}[thm]{Conjecture}

\newcommand{\Sbv}{S_{\mathrm{BV}}}
\newcommand{\Sgf}{S_{\mathrm{gf}}}
\newcommand{\sX}{{\tilde X}}
\newcommand{\sxi}{{\tilde\xi}}
\newcommand{\seta}{{\tilde\eta}}
\newcommand{\cev}[1]{{\stackrel{\leftarrow}{#1}}}
\title
{A path integral approach to the Kontsevich quantization
formula
}
\author[{A. S. Cattaneo and G. Felder}]
{ Alberto S. Cattaneo \and Giovanni Felder}
\address{A. S. C.: 
Institut f\"ur Mathematik, Universit\"at Z\"urich,
CH-8057 Z\"urich, Switzerland}
\address{G. F.:
Departement Mathematik, ETH-Zentrum, CH-8092
Z\"urich, Switzerland}
\email{asc@math.unizh.ch, felder@math.ethz.ch}
\begin{document}
\begin{abstract}
We give a quantum field theory interpretation of 
Kontsevich's deformation quantization formula for
Poisson manifolds. We show that it is given by
the perturbative expansion of the
path integral of a simple topological
bosonic open string theory. Its Batalin--Vilkovisky quantization 
yields a superconformal field theory. The associativity
of the star product, and more generally the formality
conjecture can then be understood by field theory methods.
As an application, we compute the center of the deformed algebra
in terms of the center of the Poisson algebra.
\end{abstract}\maketitle
\date{February 1999}
\section{Introduction}
In a recent paper \cite{K}, M. 
Kontsevich gave a general formula
for the deformation quantization \cite{BFFLS} of the algebra of functions
on a Poisson manifold. The deformed
product (the ``star product'') is given in terms of an expansion
reminiscent of the Feynman perturbation expansion of a two
dimensional field theory on a disc with boundary. We review 
Kontsevich's formula in Sect.\ \ref{s-2}.

The purpose of this paper is to describe this quantum field
theory explicitly. It turns out that it is a simple bosonic
topological quantum field theory on a disc $D$ with a field
$X:D\to M$ taking values in the Poisson manifold $M$ and a
one-form $\eta$ on $D$ taking values in the pull-back
$X^*(T^*M)$ of the cotangent bundle. The formula for
the star product is
\[
f\star g\,(x)=\int_{X(\infty)=x}
f(X(1))g(X(0))e^{\frac i\hbar S[X,\eta]}dX\,d\eta,
\]
where $0,1,\infty$ are three distinct points on the boundary of $D$. 
The integral is
normalized in such a way that in the case of the trivial
Poisson structure the star product
reduces to the ordinary product.
The
action $S$ is described in Sect.\ \ref{s-sigma} and was originally
studied for manifolds without boundary in
\cite{I} and  \cite{SchStr}. In particular the canonical 
quantization on the cylinder was considered.

 In the symplectic case
the above formula essentially reduces to the original Feynman path integral
formula for quantum mechanics, as pointed out to us by
H. Ooguri.

The quantization of the theory is 
somewhat subtle, due to the presence of a gauge symmetry
which only closes on shell, as already noticed in \cite{I}.
 In other words, the action $S$ is
a function of the fields
annihilated
 by a distribution of vector fields which is only integrable
on the set of critical points of $S$. As a consequence,
the BRST quantization fails and one has to resort
to the Batalin--Vilkovisky method (see for example 
\cite{BV,W1,S1,AKSZ}).

This method yields a gauge fixed action, which
turns out to have a superconformal invariance. Its
perturbative expansion around constant classical
solutions reproduces Kontsevich's formula.

As an application, we show in Sect.\ \ref{s-centre}
by quantum field theory methods that there exists a 
star product equivalent to Kontsevich's whose center consists
of the power series in $\hbar$ whose coefficients
are in the center of the Poisson algebra. 
A rigorous proof of this statement will appear elsewhere \cite{CFT}.

More generally, we may consider a path integral associated
to an arbitrary {\em poly\-vector field}, a formal sum of
skew-symmetric contravariant tensor fields of arbitrary rank,
the star product being the special case of bivector fields.
Correlation functions of boundary fields yield then a map
$U$ from poly\-vector fields to poly\-differential operators. Formal
properties of this map can be deduced from BV and factorization
methods of quantum field theory. This leads to identities,
also found by Kontsevich, which may be thought of as the
open string analog of the WDVV equations \cite{W2,DVV}. 
They may be formulated
by saying that $U$ is an $L_\infty$ morphism \cite{SchlSt, LS}. They imply
the associativity of the star product and,
in the general setting of arbitrary poly\-vector fields,
the formality conjecture \cite{K}. These constructions are explained
in Sect.\ \ref{s-Linfty}.

Although the non-rigorous quantum field theory arguments of
this paper are of course no substitute for the proofs in 
\cite{K},  this approach offers an explanation for why 
Kontsevich's construction works, and puts it in the context
of Feynman's original picture of quantization \cite{F}.
Moreover, our approach indicates the way for more general
constructions. In particular, one can consider the perturbative
expansion around a non-trivial classical solution, one can
insert a Hamiltonian and one can consider this quantum field
theory on a complex curve of higher genus.
We plan to study these variants in the future.

\medskip

\noindent{\bf Acknowledgements.} We thank H. Ooguri, whose
clarifying comments at an early stage of this work
were essential to our understanding of the problem. 
We also thank J. Fr\"ohlich and C. Schweigert for interesting
discussions and J. Stasheff for useful comments on the first
draft of this paper.
\section{The Kontsevich formula}\label{s-2}
In \cite{K}, M. Kontsevich wrote a beautiful
explicit solution to the problem of deformation quantization
of the algebra of functions on a Poisson manifold $M$. The
problem is to find a deformation of the product on the 
algebra of smooth functions on a Poisson manifold, which
to first order in Planck's constant is given by the Poisson
bracket.

If $M$ is an open set  in $\R^d$ with a Poisson structure
$\{f,g\}(x)\!=\!\sum_{i,j=1}^d
\alpha^{ij}(x)\partial_if(x)\partial_jg(x)$ given by
a skew-symmetric bivector field $\alpha$, obeying the
Jacobi identity
\begin{equation}\label{e-Jacobi}
\alpha^{il}\partial_l\alpha^{jk}+
\alpha^{jl}\partial_l\alpha^{ki}+
\alpha^{kl}\partial_l\alpha^{ij}=0,
\end{equation}
the problem is to find an associative
product $\star$ on $C^\infty(M)[[\hbar]]$, such that
for $f,g\in C^\infty(M)$,
\[
f\star g\,(x)=f(x)g(x)+\frac{i\hbar}{2}\,\{f,g\}(x)+O(\hbar^2).
\]
Kontsevich's solution\footnote{In \cite{K} $\hbar$ is what
is here $i\hbar/2$. We adopt the notation of the physics
literature and work accordingly over the complex numbers. 
With Kontsevich's conventions one may formulate the problem
over the real numbers, which in terms of the physics conventions
would mean to have an imaginary Planck constant.}
to this problem may be described
as follows. The coefficient of $(i\hbar/2)^n$ in $f\star g$ is
given by a sum of terms labeled by diagrams of order $n$. A diagram
$\Gamma$ of order $n$
is a graph consisting of $n$ vertices numbered from
$1$ to $n$ and two vertices labeled by letters
$L$ and $R$, for Left and Right. From each of the numbered vertices there
emerge two ordered oriented edges that end at numbered vertices
or at vertices labeled by letters, so that no edge starts and ends
at the same vertex. The two edges emerging from vertex $i$
are called $e^1_i,e^2_i$. They are of the form $e^a_i=(i,v_a(i))$ for
some maps $v_a:\{1,\dots,n\}\to\{1,\dots,n,L,R\}$. In fact, a 
diagram $\Gamma$ can be thought of as an ordered pair $(v_1,v_2)$
of maps $\{1,\dots,n\}\to\{1,\dots,n,L,R\}$, such that $v_a(i)$ is
never equal to $i$.

To each diagram $\Gamma$ 
of order $n$ there corresponds a bidifferential operator $D_\Gamma$
whose coefficients are differential poly\-nomials, homogeneous of
degree $n$ in the components $\alpha^{ij}$ of the Poisson structure.
The edges indicate how the partial derivatives are acting.
For instance the bidifferential operator\footnote{We use throughout the paper
the Einstein summation convention, meaning that sums over repeated indices
are understood} $(f,g)\mapsto \alpha^{ij}(x)
\partial_i f(x)\partial_jg(x)$ corresponds to the diagram with vertices $1, L, R$ and
edges $e^1_1=(1,L),e^2_1=(1,R)$. The bidifferential operator
$D_\Gamma(f\otimes g)
=\alpha^{ij}\partial_{i}\alpha^{kl}\partial_j\partial_lf\partial_kg$
corresponds to the diagram $\Gamma$ with vertices $1, 2, L, R$ and edges
$e^1_1=(1,2), e^2_1=(1,L)$, $e^1_2=(2,R), e^2_2=(2,L)$.

Kontsevich's formula is then
\[
f\star g=fg+\sum_{n=1}^\infty \left(\frac{i\hbar}2\right)^n
\sum_{\Gamma\, \mathrm{of\, order }\,n}w_\Gamma\, D_\Gamma(f\otimes g).
\]
The weight $w_\Gamma$ is the integral of a differential form
over the configuration space $C_n(H)=\{ u\in H^n, u_i\neq u_j (i\neq j)\}$
of $n$ ordered points on the
upper half plane $H$. It is defined as follows: 
for any two distinct points $z,w$ in the upper half plane
with the Poincar\'e metric $ds^2=(dx^2+dy^2)/y^2$, let
$\phi(z,w)$ be the angle between the (vertical) geodesic
connecting $z$ to $i\infty$ and the geodesic connecting 
$z$ to $w$, measured in counterclockwise direction.
Let $d\phi(z,w)=dz\frac \partial{\partial z}\phi(z,w)+
dw\frac\partial{\partial w}\phi(z,w)$ denote the differential of this angle.
Then the weight is 
\[
w_\Gamma=\frac 1{(2\pi)^{2n}n!}\int_{C_n(H)} \wedge_{i=1}^n
d\phi(u_i,u_{v_1(i)})\wedge d\phi(u_i,u_{v_2(i)}),
\]
where we set $u_L=0$ and $u_R=1$. The orientation is induced
from the product of the standard orientation of the upper half plane.

For example, we have two graphs of order one, differing in the
ordering of edges. Let us compute the weight 
of these diagrams. Let $\Gamma$ be the diagram 
with $e^1_1=(1,L), e^2_1=(1,R)$. To compute the integral over
$u=u^1+i u^2\in H$ we introduce new variables $\phi_0=\phi(u,0)$,
$\phi_1=\phi(u,1)$. As  $\mathrm{arg}(u)$ varies between $0$
and $\pi$, the angle $\phi_0$ varies from $0$ to $2\pi$. As we vary
$u$ on the half-line of constant $\phi_0$, the angle $\phi_1$ varies
between $\phi_0$ (at infinity) and $2\pi$ (at $u=0$). Thus this
change of variables is a diffeomorphism from the upper half plane
to the domain $0<\phi_0<\phi_1<2\pi$ in $\R^2$. The above description
also shows that this diffeomorphism is orientation preserving.
Thus
 $w_\Gamma=(2\pi)^{-2}\int_{0<\phi_0<\phi_1<2\pi}d\phi_0\wedge d\phi_1
=1/2$ ($\phi_j=\phi(u,j)$ and $d\phi_0\wedge d\phi_1$ is positively
oriented). The other diagram has $e^1_1=(1,R), e^2_1=(1,L)$ and gives
the same contribution with the opposite sign. Therefore the coefficient
of $(i\hbar/2)$ is 
\[
\frac12\alpha^{ij}\partial_if\partial_jg-
\frac12\alpha^{ij}\partial_jf\partial_ig
=\alpha^{ij}\partial_if\partial_jg,
\]
by the skew-symmetry of $\alpha$.

Let us conclude this section with some remarks about
involutions. The opposite
product $f\star_{\mathrm{op}}g$ is related to the product by
a change of sign of $\hbar$. Indeed, $D_\Gamma(g\otimes f)
=D_{\bar\Gamma}(f\otimes g)$ where $\bar\Gamma$ is obtained from
$\Gamma$ by exchanging $R$ and $L$, and $w_{\bar\Gamma}=(-1)^nw_\Gamma$
if $\Gamma$ is of order $n$, since $w_{\bar\Gamma}$ is the integral
of the pull-back of the differential form defining $w_\Gamma$
by the reflection about the axis $\mathrm{Re}(z)=\frac12$ which
reverses the orientation of $H$. Since the weights $w_\Gamma$ are real,
this implies that complex conjugation, extended to 
$C^\infty(M)[[\hbar]]$ by setting $\bar\hbar=\hbar$, is an antilinear
antiautomorphism for the star product.

\section{A sigma model}\label{s-sigma}
\subsection{The classical action and its symmetries}
We start by introducing a sigma model action. The perturbative
expansion of correlation functions of boundary fields on
the disc will then be related to the star product.

The model has two real bosonic fields $X,\eta$. $X$ is a map
from the disc $D=\{u\in\R^2,\, |u|\leq 1\}$ to
$M$ and $\eta$ is a differential 1-form on $D$ taking values in
the pull-back by $X$ of the cotangent bundle of $M$, 
i.e.\ a section of $X^*(T^*M)\otimes T^*D$.
In local coordinates, $X$ is given by $d$ functions $X^i(u)$
and $\eta$ by $d$ differential 1-forms $\eta_i(u)=
\eta_{i,\mu}(u)du^\mu$.

The action reads
\[
S[X,\eta]=\int_D\eta_i(u)\wedge dX^i(u)+{\frac12}\,\alpha^{ij}
(X(u))\eta_i(u)\wedge\eta_j(u).
\]
The boundary condition for $\eta$ is that for
$u\in \partial D$, $\eta_i(u)$ vanishes on vectors tangent
to $\partial D$. 

We then claim that the star product is given by the
semiclassical expansion of the
path integral\footnote{In the symplectic case, where
$\alpha$ comes from a symplectic form $\omega$, one
can integrate formally over $\eta$ and this formula
may, in the spirit of Feynman \cite{F}, be written as
\[
f\star g\,(x)=\int_{\gamma(\pm\infty)=x}
f(\gamma(1))g(\gamma(0))e^{\frac i{\hbar} \int_{\gamma}d^{-1}\omega}d\gamma\,.
\]
The integral over trajectories $\gamma:\R\to M$ is to be understood as
an expansion around the classical solution $\gamma(t)=x$, which is a
constant function of time since the Hamiltonian vanishes.
}
\[
f\star g\,(x)=\int_{X(\infty)=x}
f(X(1))g(X(0))e^{\frac i\hbar S[X,\eta]}dX\,d\eta.
\]
Here $0$, $1$, $\infty$ are any three cyclically
ordered points on the unit
circle (which we secretely view as the completed real line
by stereographic projection).
Cyclically ordered means that if we start from $0$ and move
on the circle counterclockwise we first meet $1$ and
then $\infty$. The  path integral is over all
$X:D\to M$, $\eta\in \Gamma(D,X^*(T^*M)\otimes T^*D)$ subject
to the boundary conditions $X(\infty)=x$, $\eta(u)(\xi)=0$ if
$u\in\partial D$ and $\xi$ is tangent to $\partial D$. Its
semiclassical expansion is to be understood as an expansion
around the classical solution $X(u)=x$, $\eta(u)=0$.

To evaluate this path integral we have as usual to
take gauge fixing and renormalization into
account.

This action is invariant under 
the following infinitesimal gauge transformations
with infinitesimal parameter $\beta_i$, which is
a section of  $X^*(T^*M)$ and vanishes on
the boundary of $D$:
\begin{eqnarray*}
\delta_\beta X^i&=&\alpha^{ij}(X)\beta_j,\\
\delta_\beta \eta_i&=&-d\beta_i-
\partial_i\alpha^{jk}(X)\eta_j\beta_k.
\end{eqnarray*}
This symmetry is an extension of more familiar gauge symmetries
encountered in special cases. On one extreme we have
$\alpha=0$ and the action is invariant under translations
of $\eta$ by exact one-forms on $D$. On the other
extreme we have the symplectic case where  $\alpha^{ij}$ is an
invertible matrix so that integrating formally
over $\eta$ we get the action $\int_D X^*\omega$ which
is invariant under arbitrary translations $X^i\mapsto
X^i+\xi^i$, with $\xi^i(u)=0$ on the boundary of $D$.
Another special case is the case when  $M$ is a vector
space and $\alpha$ is a {\em linear} function on $M$. In this
case $M$ is the dual space to a Lie algebra $\g$ with Kirillov--Kostant
Poisson structure. The Lie bracket of two linear functions $f,g\in\g= M^*$
is just the Poisson bracket and is again a linear function on $M$.
Then the classical action is best viewed as a function of 
a field $X$ taking values in $\g^*$ and a connection $d+\eta$
on a trivial principal bundle on $D$. 
After an integration by parts, the action becomes the ``BF action'' \cite{S2,BT}
$S=\int_D \langle X,F(\eta)\rangle$ 
where $F(\eta)$ is the curvature of $d+\eta$. In this case
the gauge transformation is the usual gauge transformation
(with gauge parameter $-\beta$)
of a connection and a field $X$ in the coadjoint representation.

In the general case, the commutator of two  gauge transformations
is a gauge transformation only on shell,
i.e., modulo the equations of motion:
\begin{eqnarray*}
{}[\delta_\beta,\delta_{\beta'}]X^i&=
&\delta_{\{\beta,\beta'\}}X^i,\\
{}[\delta_\beta,\delta_{\beta'}]\eta_i&=&
\delta_{\{\beta,\beta'\}}\eta_i
-\partial_i\partial_k\alpha^{rs}\beta_r\beta'_s
(dX^k+\alpha^{kj}(X)\eta_j).
\end{eqnarray*}
Here $\{\beta,\beta'\}_i=-\partial_i\alpha^{jk}(X)
\beta_j\beta'_k$ and $dX^k+\alpha^{kj}\eta_j=0$
is an Euler-Lagrange equation for the action $S$.
In this calculation the Jacobi identity \Ref{e-Jacobi}
plays an essential role.

Thus the gauge transformations form a Lie algebra
only when acting on critical points (classical
solutions) of $S$.

In the BRST formalism one then promotes the
infinitesimal gauge parameter $\beta_i$ to an anticommuting
ghost field (vanishing on the boundary of the disc)
and introduces the BRST operator $\delta_0$, an odd
derivation on the functions of $X,\eta,\beta$ such
that
\begin{eqnarray*}
\delta_0X^i&=&\alpha^{ij}(X)\beta_j\\
\delta_0\eta_i&=&-d\beta_i-\partial_i\alpha^{kl}(X)\eta_k\beta_l\\
\delta_0\beta_i&=&\frac12\,
\partial_i\alpha^{jk}(X)\beta_j\beta_k.
\end{eqnarray*}
Then $\delta_0$ is a differential on shell, i.e., it squares
to zero modulo the equations of motion. More precisely
we have $\delta_0^2X^i=\delta_0^2\beta_i=0$ and
$\delta_0^2\eta_i=
-\frac12\,\partial_i\partial_k\alpha^{rs}\beta_r\beta_s
(dX^k+\alpha^{kj}(X)\eta_j).
$
We assign a gradation, the ghost number, to our fields:
$\mathrm{gh}(X^i)=\mathrm{gh}(\eta_i)=0$, $\mathrm{gh}(\beta_i)=1$.
The BRST operator has then ghost number one. Additionally
we have the gradation of the fields as differential forms
on the disc, which will be denoted by deg: $\mathrm{deg}(X^i)
=\mathrm{deg}(\beta_i)=0$, $\mathrm{deg}(\eta_i)=1$.

In the case $M=\g^*$ of linear Poisson structures, the
second derivatives of $\alpha$ vanish, and the BRST operator
squares to zero.

\subsection{The Batalin--Vilkovisky action}\label{secBV}
If the BRST operator squares to zero only modulo
the equations of motion, the usual BRST
procedure to evaluate the path integral does not quite work,
since it essentially requires a well-defined cohomology
to construct physical observables. The generalization
of the BRST procedure that works in this case is the
Batalin--Vilkovisky method. The recipe is as follows.
One first adds antifields
$X^+,\eta^+,\beta^+$ with complementary ghost number and
degree as  differential forms on $D$. The
assignments of degree (from left to right) and ghost
number (from top to bottom) are given by 
\[
\begin{array}{rccc}
 & 0&1&2                \\
-2&    &    & \beta^{+i}       \\
-1&    & \eta^{+i}& X^{+}_i\\
0& X^i   &\eta_i & \\
1& \beta_i& & 
\end{array}
\]
One then looks for a Batalin--Vilkovisky action $\Sbv[\phi,\phi^+]$
of ghost number zero depending on fields $\phi^1,\phi^2,\dots$
(here $X^i,\eta_i,\beta_i$) and antifields
$\phi_1^+,\phi_2^+,\dots$, with $\mathrm{gh}(\phi^+_\alpha)=
-1-\mathrm{gh}(\phi^\alpha)$ and $\mathrm{deg}(\phi^+_\alpha)=
2-\mathrm{deg}(\phi^\alpha)$
subject to two requirements. The first requirement
is that $\Sbv[\phi,0]$ reduces to the classical action
$S[\phi]$ when the antifields are set to zero
and the second requirement is that $\Sbv$
obeys the quantum master equation 
\[
(\Sbv,\Sbv)-2i\hbar\triangle \Sbv=0.
\]
The BV Laplacian $\triangle$ and the BV antibracket are
defined as follows. 

Let us introduce temporarily a Riemannian metric on $D$, and
 denote by $\langle\ , \ \rangle_u$ the induced scalar product on
the exterior algebra of the cotangent space at $u$. The volume
form $\sqrt g du^1du^2$ will be denoted by $dv(u)$. The
Hodge star $*$ then obeys $\langle\alpha,\beta\rangle_u\,dv(u)=
\alpha\wedge*\beta$. The expression for the
Laplacian is better expressed in terms of the
Hodge dual antifields 
\[
\phi^*_\alpha=*\phi_\alpha^+.
\]
The
Laplacian of a function of fields and antifields is
\[
\triangle A=\sum_{\alpha}(-1)^{\mathrm{gh}(\alpha)}
\frac{\vec\delta^2 A}{\delta\phi^\alpha(u)\delta\phi^*_\alpha(u)}.
\]
The functional derivatives of a function of fields
and antifields, collectively denoted by $\psi^\alpha$,
are the distributions (de Rham currents)
defined
by
\[
\frac d{dt}
A(\psi+t\rho)\biggl|_{t=0}=
\int_D\langle\rho^\alpha(u),
\frac{\vec\delta A}{\delta{\psi^\alpha(u)}}\rangle_u dv(u)
=\int_D
\langle\frac {A\cev\delta}{\delta\psi^\alpha(u)}
,\rho^\alpha\rangle_u dv(u),
\]
for any test forms $\rho^\alpha$ of the same degree and ghost number
as $\psi^\alpha$.

Note that the Laplacian is the restriction of a distribution
on $D^2$ to the diagonal, and is thus a singular object in this
infinite dimensional context. It should be understood as the
limit of a suitably regularized expression.

The Laplacian obeys
\begin{equation}\label{e-LaplAB}
\triangle(AB)=\triangle(A)B+(-1)^{\mathrm{gh}(A)}(A,B)+
(-1)^{\mathrm{gh}(A)}A\triangle(B),
\end{equation}
where the {\em Batalin--Vilkovisky antibracket} is
\[
(A,B)=\sum_\alpha
\int_D\left(
\langle
\frac {A\cev\delta}{\delta\phi^\alpha(u)},
\frac{\vec\delta B}{\delta{\phi^*_\alpha(u)}}
\rangle
-
\langle
\frac {A\cev\delta}{\delta\phi^*_\alpha(u)},
\frac{\vec\delta B}{\delta{\phi^\alpha(u)}}\rangle 
\right)dv(u).
\]
This antibracket is better defined than the Laplacian,
in the sense that if $A$ and $B$ are local functionals
of the fields and antifields, such as the action $S$, then
the functional derivatives are regular distributions and
$(A,B)$ is again a local functional.
Moreover, it is independent of the choice of Riemannian metric:
it can be expressed without reference to the metric at the
cost of introducing signs:
\[
(A,B)=\sum_\alpha
\int_D \left(
\frac {A\cev\partial}{\partial\phi^\alpha}\wedge
\frac{\vec\partial B}{\partial{\phi^+_\alpha}}
-
(-1)^{\mathrm{deg}\,\phi_\alpha}
\frac {A\cev\partial}{\partial\phi^+_\alpha}\wedge
\frac{\vec\partial B}{\partial{\phi^\alpha}}\right).
\]
Here the
 derivatives of a function $A$ of fields
and antifields $\psi_\alpha$ are the distributions defined
by
\[
\frac d{dt}
A(\psi+t\rho)\biggl|_{t=0}=
\int_D\rho^\alpha\wedge
\frac{\vec\partial A}{\partial{\psi^\alpha}}
=\int_D
\frac {A\cev\partial}{\partial\psi^\alpha}
\wedge\rho^\alpha,
\]
for any test forms $\rho^\alpha$ of the same degree and ghost number
as $\psi^\alpha$.
The antibracket obeys the graded commutativity relation
\[
(A,B)=-(-1)^{(\mathrm{gh}(A)-1)(\mathrm{gh}(B)-1)}(B,A),
\]
and the Leibnitz rule
\begin{equation}\label{e-leibnitz}
(A,BC)=(A,B)C+(-1)^{(\mathrm{gh}(A)-1)\mathrm{gh}(B)}B(A,C).
\end{equation}
In the general case of field theories with non-trivial renormalization
the BV action depends on $\hbar$ through counterterms and the full
quantum master equation is solved by a recursive procedure
order by order in $\hbar$. Here, as we shall see,
 the renormalization is rather trivial and the Batalin--Vilkovisky
action satisfies
separately the equation $\triangle \Sbv=0$ and the
classical master equation
\[
(\Sbv,\Sbv)=0.
\]
 The classical master equation implies that the BV version of
the BRST operator
$\delta$ defined by $\delta A=(\Sbv,A)$ is a differential.
It obeys the Leibnitz rule
$\delta(AB)=\delta A B+(-1)^{\mathrm{gh}(A)}A\delta B$ and 
it acts on
fields and antifields by the rule
\[
\delta\phi^\alpha=(-1)^{\mathrm{gh}(\phi^\alpha)}
\frac{\vec\partial \Sbv}{\partial{\phi^+_\alpha}}
,\qquad
\delta\phi^+_\alpha=(-1)^{\mathrm{gh}(\phi^\alpha)
+\mathrm{deg}(\phi^\alpha)}
\frac{\vec\partial \Sbv}{\partial{\phi^\alpha}}.
\]
One semi-systematic 
way to find the BV action, which is under suitable
hypotheses unique up to the BV version of canonical
transformations,
is to start with the
obvious action $\Sbv^0=S+\int_D X_i^+\delta_0X^i+\eta^{+i}\wedge
\delta_0\eta_i-\beta^{+i}\delta_0\beta_i$, which has
BRST operator
$\delta=\delta_0$ and then add suitable terms, so that the
new BRST operator obeys
$\delta^2=0$. Since $\delta_0\eta^{+i}=\partial \Sbv^0
/\partial\eta_i$ contains a term proportional to  the equations of motion
(plus terms involving antifields)
which we need to cancel from $\delta_0^2\eta_i$, it is
natural to add a term quadratic in $\eta^+$ to achieve
our goal. It turns out that
\begin{eqnarray*}
\Sbv&=&\Sbv^0-\frac14\,\int_D\eta^{+i}\wedge\eta^{+j}
\partial_i\partial_j\alpha^{kl}(X)\beta_k\beta_l\\
&=&\int_D
\eta_i\wedge dX^i
+\frac12\alpha^{ij}(X)\eta_i\wedge\eta_j
+ X_i^+\alpha^{ij}(X)\beta_j
-\eta^{+i}\wedge
(d\beta_i+\partial_i\alpha^{kl}(X)\eta_k\beta_l)\\
&&
-\frac12\,\beta^{+i}
\partial_i\alpha^{jk}(X)\beta_j\beta_k
-\frac14\,\eta^{+i}\wedge\eta^{+j}
\partial_i\partial_j\alpha^{kl}(X)\beta_k\beta_l,
\end{eqnarray*}
does the job. Moreover $\Sbv$ is BRST closed (i.e., it obeys
$\delta\Sbv=0$), which is
equivalent to the classical master equation.
This is more conveniently shown in the superfield
formalism of the next subsection. 

We claim that, if the regularization is appropriate,
 $\triangle \Sbv=0$. Indeed the only terms contributing
to the Laplacian of the
BV action contain both a field and its antifield:
\begin{eqnarray*}
\triangle \Sbv
&=&
\triangle \int_D
 X_i^+\alpha^{ij}(X)\beta_j
-\eta^{+i}\wedge
\partial_i\alpha^{kl}(X)\eta_k\beta_l
-\frac12\,\beta^{+i}
\partial_i\alpha^{jk}(X)\beta_j\beta_k
\\
&=&(1-2+1)C\int_D
\partial_i\alpha^{ij}(X)\beta_j dv\\
&=&0.
\end{eqnarray*}
Here $C$ is an infinite constant. The factor takes into
account the contribution of the first term (1), of the
second term ($-2$ since the one-form $\eta_i$ has two
components)  and
the third term ($1$). In an appropriate regularization
scheme,
this cancellation is supposed to be valid before removing
the regularization, in spite of the fact that $C$ tends
to infinity.

Let us conclude this subsection by discussing the boundary
conditions of the various fields. The rule is that
Hodge dual antifields must have the same boundary conditions
as the fields. The boundary conditions for the fields are
that, for $u\in\partial D$,
$\beta_i(u)=0$ and $\eta_i(u)$ vanishes on vectors tangent to the
boundary. Thus $\beta^{+i}(u)=0$ and $\eta_i^+(u)$ vanishes on vectors
normal to the boundary.

\subsection{Superfield formalism} It turns out that
the calculations simplify if we combine our
fields and antifields into superfields. These
are functions of the even coordinates $u^1,u^2$ on $D$
and odd (anticommuting) coordinates $\theta^1,\theta^2$.
Thus a superfield $\phi$ has
the form $\phi(u,\theta)=\phi^{(0)}(u)+\theta^\mu\phi^{(1)}
_\mu(u)
+\theta^\mu\theta^\nu \frac12\phi^{(2)}_{\mu\nu}$. Its
components fields are a scalar function $\phi^{(0)}$,
a one-form $\phi^{(1)}=\phi^{(1)}_\mu du^\mu$ and a two form 
$\phi^{(2)}=\frac12\phi^{(2)}_{\mu\nu}du^\mu\wedge du^\nu$. 
The fields of total degree (degree+ghost number) zero
combine into  even superfields $\sX^i$, the
 ``super-coordinates''.
\[
\sX^i=X^i+\theta^\mu\eta_\mu^{+i}-\frac12\theta^\mu
\theta^\nu\beta_{\mu\nu}^{+i},
\]
and the fields of total degree one combine into
odd superfields $\seta_i$, the ``super-one-forms'':
\[
\seta_i=
\beta_i+\theta^\mu\eta_{i,\mu}+\frac12\theta^{\mu}
\theta^{\nu}X_{i,\mu\nu}^+.
\]
Let $D=\theta^\mu\partial/\partial u^\mu$. This
operator acts on component fields as the 
de Rham differential.
Then the BRST operator $\delta$  acts as an
odd derivation  on functions of the
superfields $\sX^i$, $\seta_i$ by the
rule
\begin{eqnarray*}
\delta\sX^i&=&D\sX^i+\alpha^{ij}(\sX)\seta_j,
\\
\delta\seta_i&=&
D\seta_i+\frac12\,\partial_i\alpha^{jk}(\sX)\seta_j\seta_k.
\end{eqnarray*}
It is easy to check that the Jacobi identity implies $\delta^2=0$.
The action of $\delta$ on components field can then
easily be evaluated by comparing coefficients
and taking into account the sign rule
 $\delta\phi=\delta\phi^{(0)}-\theta^\mu\delta\phi^{(1)}_\mu
+\frac12\theta^\mu\theta^\nu\delta\phi^{(2)}_{\mu\nu}$.
One gets
\begin{eqnarray*}
\delta X^i&=&\alpha^{ij}(X)\beta_j,
\\
\delta\eta^{+i}&=&
-dX^i-\alpha^{ij}(X)\eta_j-
\partial_k\alpha^{ij}(X)\eta^{+k}\beta_j,\\
\delta\beta^{+i}&=&
-d\eta^{+i}
-\alpha^{ij}(X)X^+_j
+\frac12\,\partial_k\partial_l\alpha^{ij}(X)
\eta^{+k}\wedge\eta^{+l}\beta_j
\\
 &&+\partial_k\alpha^{ij}(X)\eta^{+k}\wedge\eta_j
+\partial_k\alpha^{ij}(X)\beta^{+k}\beta_j.
\end{eqnarray*}
and
\begin{eqnarray*}
\delta \beta_i&=&\frac12\,\partial_i\alpha^{kl}(X)\beta_k
\beta_l,
\\
\delta\eta_i&=&
-d\beta_i-\partial_i\alpha^{kl}(X)\eta_k\beta_l-
\frac12\,\partial_i\partial_j\alpha^{kl}(X)\eta^{+j}\beta_k
\beta_l,\\
\delta X^{+}_i&=&
d\eta_{i}
+\partial_i\alpha^{kl}(X)X^+_k\beta_l
-\partial_i\partial_j\alpha^{kl}(X)
\eta^{+j}\wedge\eta_k\beta_l
+\frac12\,\partial_i\alpha^{kl}(X)\eta_k\wedge\eta_l
\\
 &&-\frac14\,\partial_i\partial_j\partial_p
\alpha^{kl}(X)\eta^{+j}\wedge\eta^{+p}\beta_k\beta_l
-\frac12\,\partial_i\partial_j\alpha^{kl}(X)\beta^{+j}\beta_k
\beta_l.
\end{eqnarray*}
This BRST operator coincides with the one obtained by the
Batalin--Vilkovisky procedure. 
The Batalin--Vilkovisky action is the integral \[S_{BV}=\int_D L^{(2)}.\]
of the two-form part $L^{(2)}=\int d^2\theta L$ of 
\[
L=\seta_iD\sX^i+\frac12\,\alpha^{ij}(\sX)\seta_i
\seta_j.
\]
It is BRST closed, i.e., it obeys the master equation. In fact one has
\[
\delta L= D(\seta_iD\sX^i),
\]
so that $\delta L^{(2)}$ is the differential of a one form which 
vanishes along the boundary.
\subsection{The gauge fixed action}
We compute the path integral in the Lorentz-type gauge 
$d{*\eta_i}=0$. The Hodge $*$ operator (alias the almost complex structure) 
depends on the conformal structure
and the orientation of $D\subset \R^2$: in terms of the standard coordinates,
$*du^1=du^2$, $*du^2=-du^1$. 

Let us shortly recall the main idea of 
the Batalin--Vilkovisky formalism in the general setting of
\ref{secBV}. For any function $\Psi$, the ``gauge fixing fermion'',
of the fields of
ghost number $-1$, one considers the 
integral $\int_L \mathcal{O} e^{\frac i\hbar \Sbv}$ for an observable
$\mathcal{O}$, i.e., a function of fields and antifields which is
closed with respect to the quantum BRST operator $\Omega=-i\hbar\triangle +\delta$:
\[
\Omega\,\mathcal{O}=-i\hbar \triangle\mathcal{O}+(\Sbv,\mathcal{O})=0.
\]
The integral is   taken over
the ``Lagrangian'' submanifold $L$ defined by the
equations $\phi_\alpha^+=\vec\partial_{\phi^\alpha}\Psi$. 
Using formally the master equation and the fact that $\mathcal{O}$ is
BRST closed, one then see that
 these integrals are invariant under variations of $\Psi$ and
thus ``equal'' to the original (ill-defined) path integral with
action $S[\phi]=\Sbv[\phi,0]$, which is what one gets if $\Psi=0$.

The problem is then to find a function $\Psi$ which makes
the integral well-defined, at least as a perturbative series.
One way to do this is to add new fields, called antighosts
and Lagrange multipliers together with their antifields, 
and choose $\Psi$ as the scalar
product of the antighost and the gauge fixing condition. The 
action for these new fields is the simplest and is
added to the Batalin--Vilkovisky action $\Sbv$.

Let us do this in the case at hand. We
first introduce new anticommuting scalar fields (antighosts) $\gamma^i$ 
of ghost number $-1$ on $D$, and scalar
Lagrange multiplier fields $\lambda^i$ of ghost number zero,
together with their antifields $\gamma^+_i$, $\lambda^+_i$. The
boundary condition for $\lambda^i$ is Dirichlet: $\lambda^i(u)=0,
u\in\partial D$, and $\gamma^i$ is constant on the boundary. The action
for these fields and antifields is $-\int_D\lambda^i\gamma^+_i$ and 
is just added to the BV action. The BRST operator acts then as
\[
\delta\lambda=\delta\gamma^+=0,\qquad \delta\lambda^+=-\gamma^+,
\qquad \delta\gamma=\lambda.
\]
Clearly the new action also obeys the
master equation.
The gauge fixing condition $d{*\eta}=0$ is encoded in the
gauge fixing fermion
$\Psi=-\int_Dd\gamma^i{{*\eta_i}}$.
On the Lagrangian submanifold we then have $X^+=\beta^+=\lambda^+=0$,
$\gamma^+_i=d{*\eta_i}$ plus a boundary term whose form will not
matter,
 and $\eta^{+i}=*d\gamma^i$. The boundary condition for $\gamma^i$ 
was chosen so as to fulfill the boundary condition for $\eta^+$ (vanishing
on normal vectors).
The gauge fixed action is then
\begin{eqnarray*}
\Sgf&=&\int_D\eta_i\wedge dX^i+{\frac12}\,\alpha^{ij}
(X)\eta_i\wedge\eta_j-*d\gamma^i\wedge(d\beta_i+
\partial_i\alpha^{kl}(X)\eta_k\beta_l)
\\
&&-\frac14{*d\gamma}^i\wedge {*d\gamma}^j
\partial_i\partial_j\alpha^{kl}(X)\beta_k\beta_l-\lambda^id{*\eta_i}.
\end{eqnarray*}

\subsection{Superconformal invariance of the gauge fixed action} 
The original action is invariant under arbitrary diffeomorphisms
of the disc. As the gauge fixing condition depends on a choice of
conformal structure, the gauge fixed action is only invariant
under conformal diffeomorphisms. In fact this invariance is
part of a (twisted) superconformal invariance, as we now show.
For each vector field $\epsilon(u)=\epsilon^\mu(u)
\frac\partial{\partial u^\mu}$ on $D$, tangent to the boundary
on $\partial D$, we introduce  an odd derivation
$\bar\delta_\epsilon$, depending linearly on $\epsilon$,
 on functions of our fields:
\[
\begin{array}{lll}
\bar\delta_\epsilon X^i=i(\epsilon){*d\gamma^i},&
\bar\delta_\epsilon\lambda^i=-i(\epsilon)d\gamma^i,
\\ && \\
\bar\delta_\epsilon\beta_i=i(\epsilon)\eta_{i}, &
\bar\delta_\epsilon\eta_{i}=0, &
\bar\delta_\epsilon\gamma^i=0.
\end{array}\]
Here $i(\epsilon)$ is the interior multiplication of a differential
form on $D$ with a vector field $\epsilon$.
A straightforward calculation shows that these
derivations, together with the BRST operator obey the
twisted supersymmetry algebra relations 
\[
{}[\delta,\delta]_+=[\bar\delta_\epsilon,\bar\delta_\epsilon']_+=0,
\qquad [\delta,\bar\delta_\epsilon]_+=-\mathcal{L}_\epsilon,
\]
modulo the equations of motion for $\Sgf$,
with Lie derivative $\mathcal{L}_\epsilon=i(\epsilon)\circ d
+d\circ i(\epsilon)$.

The gauge fixed action obeys $\delta \Sgf=0$ and 
\[
\bar\delta_\epsilon\Sgf
=\int_D\eta_i\wedge(\mathcal{L}_\epsilon{*d\gamma^i}-*\mathcal{L}_\epsilon
d\gamma^i).
\]
The latter expression vanishes if $\epsilon$ is conformal,
 i.e., if $\mathcal{L}_\epsilon$ commutes with $*$. Conformal
vector fields on the disc form a three dimensional Lie algebra,
isomorphic to $\mathfrak{su}(1,1)$.

\subsection{BRST cohomology classes}
Observables can be
obtained from differential forms on $M$. To
a differential $p$-form $\omega=\omega_{i_1,\dots,i_p}(x)
dx^{i_1}\wedge\cdots\wedge dx^{i_p}$ and a point
$u$ on the boundary of the disc $D$, one associates
the observable
\[
\hat\omega(u)
=\omega_{i_1,\dots,i_p}(X(u))\gamma^{i_1}(u)\cdots\gamma^{i_p}(u).
\]
In general, functions of the components of $\sX(u)$ with $u$ on the
boundary are observables, since with our boundary conditions we
have $\delta\sX=0$ on the boundary and the Laplacian of a function
depending only a field but not on its antifield or on an antifield
but not its field, vanishes. More general observables are considered
in Sect.\ \ref{s-centre}.

The gauge fixed action still has a (finite dimensional)
residual infinitesimal symmetry $X^i\mapsto X^i+a^i$, $\gamma^i
\mapsto \gamma^i+g^i$ where $a_i$, $g_i$ are constant 
functions on the disc. This translates into zero modes
in the integration over the fermions $\gamma^i$, and it follows
that the only observables that have non-zero integral 
have ghost number $-\mathrm{dim}(M)$.

\subsection{Feynman rules}
The Feynman perturbation expansion in powers of $\hbar$ around the
classical solution $X(u)=x$, $\eta(u)=0$ can be now computed.  Thus we
write $X(u)=x+\xi(u)$ with a fluctuation field $\xi(u)$ with
$\xi(\infty)=0$. The Feynman propagators can then be deduced from the
``kinetic'' part 
\begin{eqnarray*}
\Sgf^{0}&=&\int_D\eta_i\wedge d\xi^i
-*d\gamma^i\wedge
d\beta_i-\lambda^id{*\eta_i}
\\
&=&
\int_D\eta_i\wedge (d\xi^i+*d\lambda^i)+\beta_id{*d}\gamma^i.
\end{eqnarray*} 
of the gauge fixed action. The other
terms of $\Sgf$ are considered as perturbations. Thus we have
to invert the operators $d\oplus*d:\Omega^0(D)\oplus \Omega_0^0(D)
\to\Omega^1(D)$ and $d{*d}:\Omega^0(D)\to\Omega^2(D)$. Here
$\Omega^p(D)$ is the space of smooth $p$-forms on $D$ and $\Omega^0_0(D)$
denotes the space of functions with Dirichlet boundary conditions
$\lambda^i(u)=0, u\in\partial D$. Both operators are surjective but
have a one-dimensional kernel consisting of constant functions.
Inverses (modulo these kernels) are  
integral operators: to describe them it is useful to
map conformally the disc onto the upper half plane $H_+$ and
use the standard complex coordinate of $H_+$. 
The integration kernel of
$(d{*d})^{-1}$ is the Green function 
$\frac1{2\pi}\psi(z,w)$, with
\[
\psi(z,w)=\ln\left|\frac{z-w}{z-\bar w}\right|.
\]
The integration kernel of $(d\oplus {*d})^{-1}$ is the Green function
$G(w,z)=\frac1{2\pi}(*d_z\psi(z,w)\oplus d_z\phi(z,w))$, where
$d_z=dz\frac\partial{\partial z}+d\bar z\frac{\partial}{\partial \bar z}$ is
the differential with respect to $z$ and
\[
\phi(z,w)=\frac1{2i}\ln\frac{(z-w)(z-\bar w)}{(\bar z-\bar w)(\bar z-w)},
\]
We have $d_w{*d_w}\psi(z,w)=d_w{*d_w}\phi(z,w)=2\pi\delta_z(w)$ where
$\delta_z(w)$ is the Dirac distribution two-form, and the boundary
conditions for $w\in\partial H_+$ are Dirichlet for $\psi$ and
Neumann for $\phi$.
The propagators are then
\[
\langle\gamma^k(w)\beta_j(z)\rangle=
\frac{i\hbar}{2\pi}\delta_{j}^k\psi(z,w),
\qquad
\langle\xi^k(w)\eta_j(z)\rangle=
\frac{i\hbar} {2\pi}\delta^{k}_jd_z\phi(z,w),
\]\[
\langle\lambda^k(w)\eta_j(z)\rangle
=\frac{i\hbar}{2\pi} \delta^{k}_j{*d_z}\psi(z,w).
\]
 Note that $*d_w\psi(z,w)=d_w\phi(z,w)$ so that
$\langle{*d}\gamma^k(w)\beta_j(z)\rangle=
\delta_j^k\frac{i\hbar}{2\pi }
d_w\phi(z,w)$. 
It follows that the  propagators combine into a {\em superpropagator}
\[
\langle\xi^k(w)\eta_j(z)\rangle+\langle{*d}\gamma^k(w)\beta_j(z)\rangle
=\frac{i\hbar}{2\pi}\delta^{k}_jd\phi(z,w),
\]
where $d=d_z+d_w$. In terms of superfields
$\seta_j(z,\theta)=\beta_j(z)+\theta^\mu\eta_{j,\mu}(w)$,
$\sxi^k(w,\zeta)=\xi^k(w)+\zeta^\mu\eta^{+j}_\mu(w)$, with
$\eta^{+j}=*d\gamma^j$, the superpropagator is
\[\langle\sxi^k(w,\zeta)\seta_j(z,\theta)\rangle
=\frac{i\hbar}{2\pi}\delta^{k}_jD\phi(z,w),
\]
where $D=\theta^\mu\frac\partial{\partial z^\mu}
+\zeta^\mu\frac\partial{\partial w^\mu}$.

The perturbation expansion is then obtained by writing
$\Sgf=\Sgf^{0}+\Sgf^{1}$ and expanding:
\[
\int e^{\frac i\hbar\Sgf}\mathcal O=\sum_{n=0}^\infty\frac{i^n}{
\hbar^nn!}
\int e^{\frac i\hbar\Sgf^{0}}(\Sgf^{1})^n\mathcal O.
\]
This expression is calculated using the Wick theorem
for Gaussian integrals
\begin{eqnarray*}
\lefteqn{\int e^{\frac i\hbar\Sgf^{0}}
\sxi^{k_1}(w_1,\zeta_1)\cdots\sxi^{k_N}(w_N,\zeta_N)
\seta_{j_1}(z_1,\theta_1)\cdots\seta_{j_N}(z_N,\theta_N)
\delta_x(X(\infty))}\\
&=&
\sum_{\sigma\in S_N}
\langle\sxi^{k_{\sigma(1)}}(w_{\sigma(1)},\zeta_{\sigma(1)})\seta_{j_1}(z_1,
\theta_1)\rangle
\cdots
\langle\sxi^{k_{\sigma(N)}}(w_{\sigma(N)},\zeta_{\sigma(N)})\seta_{j_N}(z_N,
\theta_N)\rangle.
\end{eqnarray*}
The normalization of the integral is such that
$\int\exp(\frac i\hbar\Sgf^{0})\delta_x(X(\infty))=1$, 
so that for $\alpha=0$ the
star product coincides with the ordinary product.
Here $\delta_x(X(t))=\prod_{i=1}^d\delta(X^i(t)-x^i)\gamma^i(t)$ fixes the
value of the zero modes (constant functions) of $X$ and 
the $\gamma$'s are needed since the integral is otherwise zero,
owing to the presence of zero modes in the integration over $\gamma$.
More generally we could insert instead of the delta distribution a
factor $\rho(X(\infty))\gamma^1(\infty)\cdots\gamma^d(\infty)$,
for some top differential form $\omega=\rho(x)dx^1\cdots dx^d$
on $M$, resulting in a factor $\int_M\omega$ in the right-hand side.
The Feynman perturbation expansion is then obtained by expanding
the interaction term $\Sgf^1$ and the observable in powers of $\sxi$
$\seta$. This gives the vertices
\begin{equation}\label{e-vert}
\Sgf^1=\frac12\int_D\int d^2\theta\sum
_{k=0}^\infty\frac1{k!}
\partial_{j_1}\cdots\partial_{j_k}\alpha^{ij}(x)
\sxi^{j_1}\cdots\sxi^{j_k}\seta_{i}\seta_{j}.
\end{equation}
Here the Berezin integral selects the two-form part of
the superfield.
We consider the observable of the correct ghost number
\begin{equation}\label{e-O}
\mathcal{O}=f(\sX(1))g(\sX(0))\delta_x(X(\infty)),
\end{equation}
where $f,g\in C^\infty(M)$ . Then
expanding $f$ and $g$ in powers of $\sxi$ we get
an expansion in Feynman diagrams\footnote{Actually, the arguments
of $f$ and $g$ of our original integral are $X(1)$, $X(0)$,
 rather
than the superfields. However, the additional terms in $\sX(1),
\sX(0)$
do not 
contribute to the integral since their are of negative ghost
 number}.
 The terms with $n$
vertices \Ref{e-vert} are then labeled by the Kontsevich diagrams
$\Gamma$ of order $n$, but possibly with tadpoles, i.e.,
lines that start and
end at the same vertex. 
The term labeled by a diagram $\Gamma$
with lines $(j,v_1(j))$, $(j,v_2(j))$, $j=1,\dots,n$ is
$D_\Gamma(g\otimes f)$ (see Sect.\ \ref{s-2}) times
\[
\frac1{n!}
\left(
\frac i \hbar
\right)^n
\frac1{2^n}
\left(\frac {i\hbar}{2\pi}\right)^{2n}
\int \wedge_{j=1}^n
d\phi(u_j,u_{v_1(j)})\wedge d\phi(u_j,u_{v_2(j)})=
(-1)^nw_\Gamma,
\]
The factor $1/\prod k_j!$, where $k_j$ is the number
of lines pointing to $j$, is compensated by the fact that
there are as many terms in the Wick theorem which give the
same contribution because $k_j$ arguments of $\sxi$ are equal
to each other.

As explained at the end of Sect.\ \ref{s-2}, we have
$(-1)^nw_\Gamma D_\Gamma(g\otimes f)=w_{\bar\Gamma}
D_{\bar\Gamma}(f\otimes g)$,
where $\bar\Gamma$ is $\Gamma$ with $R$ and $L$ interchanged.
Thus the product obtained here
 coincides with Kontsevich's, except that it also involves
tadpoles diagrams. These have to be considered separately,
and require (finite) renormalization, which we proceed to discuss.

\subsection{Renormalization}\label{ss-ren}
In the perturbation expansion described above, all integrals
are absolutely convergent except for those containing tadpole
diagrams, which are diagrams with one edge connecting a vertex
to itself. The corresponding amplitude contains an ill-defined
factor $d\phi(z,z)$, the superpropagator taken at coinciding points.

To make sense of this expression we introduce a point-splitting
regularization
and define $d\phi(z,z)$ as the limit
\[
d\phi(z,z)=\kappa(z;\zeta)
=\lim_{\epsilon\to 0}d\phi(z,z+\epsilon\zeta(z)).
\]
Here $\zeta(z)$ is a vector field on $D$ which does not vanish
in the interior of $D$. This limit exists but depends on
the regularizing vector field $\zeta(z)$. Indeed, if we
write $\zeta(z)=r(z)e^{i\vartheta(z)}$ in polar coordinates,
then
\[
\kappa(z;\zeta)={d\vartheta(z)}.
\]
Thus the Feynman integrals have a finite renormalization ambiguity.
One way to fix it is to add a counterterm 
\begin{equation}\label{e-ct}
S_{\mathrm{c.t.}}=\frac{i\hbar}{2\pi}\int_D\int d^2\theta\,
\partial_i\alpha^{ij}(\sX)
\seta_j\tilde\kappa,\qquad
\tilde\kappa=\theta^\mu\kappa_\mu,
\end{equation}
(or more simply choose the slightly singular
$\vartheta=$ constant) which removes the tadpoles diagrams,
and one gets precisely the Kontsevich formula. One easily checks
that the action with the addition of the counterterm still obeys
the classical master equation and, by the same argument as at the
end of \ref{secBV}, also the quantum master equation.

%%%%%%%%%%%
\section{Central functions}\label{s-centre}
Using a non-rigourous quantum field theory argument based on BRST cohomology,
we can prove the following claim:
\begin{quote}
{\it There is a star product, equivalent to Kontsevich's, so
that every function that is central in the Poisson algebra 
is also central for the star product.}
\end{quote}
Two star products $\star$, $\star^\prime$ corresponding to the same Poisson bracket
are called equivalent if
there is a series $R=R_0+\hbar R_1+\hbar^2 R_2+\cdots$ with $R_i$ 
differential operators and $R_0=\mathrm{Id}$, such that
$f\star^\prime g=R^{-1}(Rf\star Rg)$.
%****
The argument will also give us a formula  for $R$, see \Ref{e-R}.
%****

Observe first that the BRST variation of a function on $M$ is given by
\[
\delta f(\sX) = \partial_if(\sX)\,(D\sX^i+\alpha^{ij}(\sX)\seta_j)=
Df(\sX) + \seta_j \alpha^{ij}(\sX) \partial_if(\sX).
\]
If $f$ is central in the Poisson algebra---that is, if $\alpha^{ij}(X)
\partial_j f(X)=0$---then the second term on the right hand side vanishes.\footnote
{In presence of the counterterm \Ref{e-ct} the BRST operator
is modified and we get an additional term $(i\hbar/2\pi)
\tilde
\kappa\partial_j(\alpha^{ji}\partial_if)(\sX)$ in the formula for
$\delta f$. This term also vanishes if $f$ is central in the
Poisson algebra.} 
Writing $f(\sX)$ in components, $f(\sX) = f(X) + \theta^\mu\eta_\mu^{+i}
\partial_i f(X)+\dotsb$, 
we get the descent equations
\begin{align*}
\delta f(X) &= 0,\\
\delta(\eta^{+i} \partial_i f(X)) &= -df(X).
\end{align*}
The first equation means that $f(X)$ is an observable. Therefore,
the expectation value
\begin{equation}\label{e-99}
h(u;f,g)(x) = \int
f(X(u))g(X(0))\delta_x(X(\infty))e^{\frac i\hbar S}
\end{equation}
is well defined for any $u$ in the upper half plane. Observe that we put no
additional hypotheses on $g$, so that---as in the previous
sections---$g(X(v))$ is an observable only if $v$ is on the real axis.

The second descent equation may then be used to prove that $h$ is independent
of $u$. In fact, denoting by $d$ the exterior derivative on the upper half
plane, we get
\begin{multline*}
dh(u;f,g)(x) = \int
df(X(u))g(X(0))\delta_x(X(\infty))e^{\frac i\hbar S}=\\=
-\int
\delta[\eta^{+i}(u) \partial_i f(X(u))]
g(X(0))\delta_x(X(\infty))e^{\frac i\hbar S}=\\=
-\int
\delta[\eta^{+i}(u) \partial_i f(X(u))
g(X(0))\delta_x(X(\infty))]e^{\frac i\hbar S}=0.
\end{multline*}
Observe that to obtain the last equality one must also check that
\[\triangle[\eta^{+i}(u) \partial_i f(X(u))g(X(0))\delta_x(X(\infty))]=0.\]
As a consequence, we get eventually
%****
\[
Rf\star g\,(x)=\lim_{\epsilon\downarrow0}
h(1+i\epsilon;f,g)(x)=
\lim_{\epsilon\downarrow0}h(-1+i\epsilon;f,g)(x)=g\star Rf\,(x),
\]
Here for $v$ on the real axis, 
$Rf(X(v))=\lim_{\epsilon\downarrow0}f(X(v+i\epsilon))$ is
the limit of the observable $f(X(u))$ defined on the upper half plane
as $u$ tends to the real axis. We claim that $Rf$ is given by the one-point
function
\begin{equation}\label{e-R}
Rf(x)=\int
f(X(u))\delta_x(X(\infty))e^{\frac i\hbar S}=f(x)+O(\hbar^2)
\end{equation}
for any point $u$ not on the boundary. This is based on the following
factorization argument: if in the integral $h(u;f,g)$ the point
$u$ approaches the boundary, it is as if we considered the integral
on two discs connected by a small bridge, with $u$ in the middle of
one disc and the insertion point for the observable $g$ on the boundary
of the other. In the limit one obtains a path integral for a disc with
two points (and the point at infinity) on the boundary. One
point is the insertion point for $g$ and at the other
the result of the path integral on the disc with one point in
the interior is inserted. See Fig.~1 for the case when $u$ approaches
$-1$.
\begin{figure}[h]
{\scalebox{.5}{\includegraphics{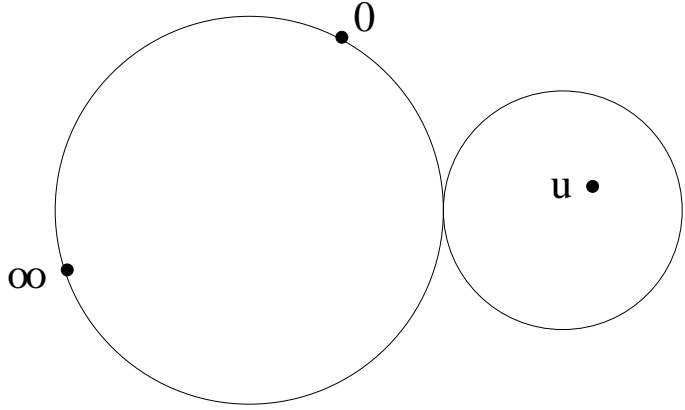}}}
\caption{The expectation value \Ref{e-99} in the limit as  $u$ approaches the
boundary reduces to a path integral on this surface}
\end{figure}
This argument can be made precise looking at the perturbation
expansion \cite{CFT} with the result that
$Rf=f+\hbar^2R_2f+\cdots$, with $R_j$ differential operators.
%****

Thus $f$ is central for the star product 
$g\star^\prime h=R^{-1}(Rg\star Rh)$,  proving
 the claim at the beginning of the section.

%%%%%%%%%%%%
{}Using this result we may strengthen our claim:
\begin{quote}
{\it The center of $C^\infty(M)[[\hbar]]$ with the star product $\star^\prime$ 
is 
$Z[[\hbar]]$ where $Z=\{f\in C^\infty(M)\,|\,\{f,\cdot\}=0\}$ is
the center of the Poisson algebra $C^\infty(M)$.}
\end{quote}
The proof goes as follows. We need to show that if $f=f_0+\hbar f_1+\cdots$
is central for $\star^\prime$ then all coefficients $f_i$ are in $Z$.
If $\{f,g\}=0$ for arbitrary $g\in C^\infty(M)$, then in particular the
coefficient of $\hbar$ vanishes, so $\{f_0,g\}=0$ and $f_0\in Z$. But this
implies, by what we showed above, that $f_0$ is central for the star product.
Thus $(f-f_0)/\hbar=f_1+\hbar f_2+\cdots$ is central. By proceeding in this way we
see that $f_1$, $f_2, \dots$ are all in $Z$.

\section{$L_\infty$ morphism and formality} \label{s-Linfty}

\subsection{The general path integral as a map from poly\-vector fields
to poly\-differential operators}
The path integral we considered so far is a special case of
the following general construction. A {\em poly\-vector field}
of degree $p$ is a section of $\wedge^{p+1}TM$, i.e., a
skew-symmetric contravariant tensor field of rank $p+1$,
\[\frac1{(p+1)!}\alpha^{j_0,\dots,j_p}(x)
\frac{\partial}{\partial x^{j_0}}
\wedge
\cdots
\wedge
\frac{\partial}{\partial x^{j_p}}.\]
 A {\em poly\-vector
field} is a sum $\alpha=\sum_{p=0}^{d-1}\alpha^{(p)}$ of 
poly\-vector fields of all nonnegative degrees. The space of poly\-vector fields
is denoted by $T_{\mathrm{poly}}(M)$.

For a multi-index $I=(i_1,\dots,i_d)\in\Z_{\geq 0}^{d}$,
let $\partial_I=\prod_k (\partial/\partial x^k)^{i_k}$.
A {\em poly\-differential operator} of degree $m$ is an operator
 $V:C^\infty(M)^{\otimes m+1}\to C^\infty(M)$,
of the form $V(f_0\otimes\cdots\otimes f_m)(x)
=\sum V_{I_0,\dots, I_m}(x)\partial_{I_0}f_0
\cdots \partial_{I_m}f_m$, with a finite sum over  sets
of multi-indices $I_j$. A poly\-differential operator is
a formal sum of poly\-differential operators of arbitrary nonnegative
degrees.
The space of poly\-differential operator is denoted by
$D_{\mathrm{poly}}(M)$.

To a poly\-vector field $\alpha$ we may associate a function
of fields and 
antifields: 
\[
S=S_0+S_\alpha\, ,
\]
with
\[
S_0=\int_D\int d^2\theta\,\seta_j D\sX^j-\int_D\lambda^i\gamma^+_i,
\]
as above, and
\[
S_\alpha=
\sum_{p=0}^{d-1}\int_D\int d^2\theta\,
\frac1{(p+1)!}\alpha^{j_0,\dots,j_p}(\sX(u,\theta))
\seta_{j_0}(u,\theta)\cdots\seta_{j_p}(u,\theta).
\]
We may then consider
 correlation functions of boundary fields associated
to the functions $f_0,\dots,f_m$ on $M$.
\begin{eqnarray*}
U(\alpha)(f_0\otimes\cdots\otimes f_m)(x)
&=&
\int e^{\frac i\hbar(S_0+S_\alpha)}\mathcal{O}_x(f_0,\dots,f_m),
\\
\mathcal{O}_x(f_0,\dots,f_m)&=&\int_{B_m}
 [f_0(\sX(t_0,\theta_0))
\cdots f_m(\sX(t_m,\theta_m))]^{(m-1)}\delta_x(X(\infty)).
\end{eqnarray*}
The path integral is, as before, the integral over the Lagrangian
submanifold in the space of fields and antifields determined
by our gauge condition $d{*\eta}=0$. The integral
over the $t_i$ is the integral over the $m-1$ form part (the coefficient of
$\theta_1\cdots\theta_{m-1}$) of
the integrand over the simplex
$1=t_0>t_1>\dots>t_{m-1}>t_m=0$, with the orientation given by the
volume form $dt_1\wedge\cdots\wedge dt_{m-1}$.
 It may be viewed as an integral over the moduli
space $B_m$ 
of $m+1$  cyclically ordered points on the circle modulo conformal 
transformations.
More explicitly
\[
\mathcal{O}_x(f_0,\dots,f_m)=
\int_{1>t_1>\dots>t_{m-1}>0}f_0(X(1))\,\prod_{k=1}^{m-1}
\partial_{i_k}f(X(t_k))\eta^{+i_k}(t_k)\,f_m(X(0))\delta_x(X(\infty)).
\]
Expanding the path integral in powers of $\hbar$ as in the previous section, we get a
map $U$ that associates, to each poly\-vector field $\alpha$,  a formal
power series whose coefficients are 
 poly\-differential 
operators.

The perturbative expansion has the form
$U(\alpha)=\sum_{n=0}^\infty 
U_n(\alpha,\dots,\alpha;\hbar)$. Here $U_n$ is a multilinear
function of $n$ arguments in $T_{\mathrm{poly}}(M)$
 with values $D_\mathrm{poly}(M)$. 
The formula for $U_n$ is
\[
U_n(\alpha_1,\dots,\alpha_n;\hbar)(f_0\otimes\cdots\otimes f_m)(x)
=
\int e^{\frac i\hbar S_0}
\frac i\hbar S_{\alpha_1}\cdots \frac i\hbar S_{\alpha_n}
\mathcal{O}_x(f_0,\dots,f_m).
\]
Suppose now that, for $i=1,\dots, n$,
 $\alpha_i$ is homogeneous of degree $p_i$. Then $S_{\alpha_i}$ is
the integral of the two-form component of $\frac1{(p_i+1)!}
\alpha_i^{j_0\dots j_{p_i}}
(\sX)\seta_{j_0}\dots\seta_{j_{p_i}}$, and has thus ghost number
$p_i-1$. This has two consequences: first, since the
integral over the $t_i$ picks the $m-1$ form component of
$\prod f_i(\sX(t_i))$, which has ghost number $1-m$, we
have the ghost number condition
$
1-m+\sum_{i=1}^n(p_i-1)=0$ or
\begin{equation}
\label{e-gh} 
m=1-n+\sum_{i=1}^n p_i,
\end{equation}
for the path integral to be non-zero. This means that
$U_n$ is a map of degree $1-n$ from $T_{\mathrm{poly}}(M)^{\otimes n}$
to $D_{\mathrm{poly}}(M)$.
Using this formula we may compute the dependence on $\hbar$ of
 this integral: 
the path integral has an overall $\hbar$ to the
power $-n+\sum (p_i+1)=n+m-1$ (each vertex has $1/\hbar$ and
each propagator has an $\hbar$), and we have
\[
U_n(\alpha_1,\dots,\alpha_n;\hbar)(\otimes_0^mf_i)=(i\hbar)^{n+m-1}
U_n(\alpha_1,\dots,\alpha_n)(\otimes_0^mf_i),
\]
with $U_n(\alpha_1,\dots,\alpha_n)=U_n(\alpha_1,\dots,\alpha_n;\hbar=1/i)$
independent of $\hbar$. 

The second consequence is that 
\begin{equation}\label{e-sym}
U_n(\dots,\alpha_i,\dots,\alpha_j,\dots)
=(-1)^{(p_i-1)(p_j-1)}U_n(\dots,\alpha_j,\dots,\alpha_i,\dots),
\end{equation}
i.e., $U_n$ is symmetric in a graded sense.

\subsection{Special cases}
Let us consider in detail some special cases. For $n=0$, $U_n$
is a poly\-differential operator of degree $m=1$, and 
\[U_0(f_0\otimes f_1)=\int e^{\frac i\hbar S_0}f_0(\sX(1))
f_1(\sX(0))\delta_x(X(\infty))=f_0(x)f_1(x)
\]
 is the
undeformed product on $C^\infty(M)$.

If $n=1$ and $\alpha$ is a poly\-vector field of degree $p$ then
$U_1(\alpha)$ is of degree $p$. Let
\[
\alpha=\frac1{(p+1)!}\alpha^{j_0,\dots,j_p}(x)
\frac{\partial}{\partial x^{j_0}}
\wedge
\cdots
\wedge
\frac{\partial}{\partial x^{j_p}},
\] 
with $\alpha^{j_0,\dots,j_p}$ antisymmetric.
The Wick theorem yields in this
case
\[
U_1(\alpha;\hbar)(f_0\otimes\cdots\otimes f_p)(x)
=\frac i\hbar\left(\frac{i\hbar}{2\pi}\right)^{p+1}
\,I_p\,
 \alpha^{j_0,\dots,j_p}\partial_{j_0}f_0(x)
\cdots
\partial_{j_p}f_p(x).
\]
Here $I_p$ is the integral
\[
I_p=\int d\phi(u,1)\wedge d\phi(u,t_1)\wedge\cdots
\wedge d\phi(u,0),
\]
over $u=u^1+iu^2\in H$ and $1>t_1>\cdots>t_{p-1}>0$,
with orientation given by the form $du^1\wedge du^2\wedge
dt_1\dots dt_{p-1}$. To compute this integral we proceed as
in Sect.\ \ref{s-2} and introduce new
variables $\phi_0=\phi(u,1)$, $\phi_j=\phi(u,t_j)$ ($j=1,\dots, m-1$)
and $\phi_p=\phi(u,0)$. In the new variables the integration is
over the region $2\pi>\phi_0>\cdots>\phi_p>0$. We claim that the
Jacobian of the change of variables is $(-1)^p$. This follows from
the fact that $d\phi(u,0)\wedge d\phi(u,1)=Jdu^1\wedge du^2$ with
$J>0$ and that $\partial\phi(u,t)/\partial t>0$. 
Hence,
\begin{eqnarray*}
\int d\phi(u,1)\wedge d\phi(u,t_1)\wedge\cdots
\wedge d\phi(u,0)&=&(-1)^p\int_{2\pi>\phi_0>\dots>\phi_p>0}
d\phi_0\cdots d\phi_p
\\
&=&\frac{(-1)^p(2\pi)^{p+1}}{(p+1)!}\, ,
\end{eqnarray*}
and we obtain
\[
U_1(\alpha)(f_0\otimes \cdots\otimes f_p)(x)=
\frac{(-1)^{p+1}}{(p+1)!}\alpha^{j_0,\dots,j_p}(x)
\partial_{j_0}f_0(x)\cdots\partial_{j_p}f_p(x).
\]
\subsection{$U$ is an $L_\infty$ morphism}
The formal properties of the map $U$ can be deduced using
the main trick of the BV formalism, which is to use the
fact that
the integral of the Laplacian of anything is zero. In our
situation we have, with $S_0=\int_D\int d^2\theta
 \seta_iD\sX^i-\int_D\lambda^i\gamma^+_i$, and
$\alpha_j$ ($j=1,\dots,n$) 
homogeneous polyvector fields of degree $p_j$,
\[
\int\triangle \left( e^{\frac i\hbar S_0}\prod_{i=1}^n S_{\alpha_i}
\mathcal{O}_x(f_0,\dots,f_m)\right)=0.
\]
To evaluate the left-hand side we use \Ref{e-LaplAB}
 and $\triangle S_0=\triangle S_\alpha=0$. Also, we
use the fact that $(S_0,S_\alpha)$ is proportional to
$\int_D\int d^2\theta D(\alpha^{j_0\dots j_p}(\sX)\seta_{j_0}\cdots\seta_{j_p})$
which vanishes because of the boundary conditions for $\seta_j$.
Thus 
we get
\begin{eqnarray*}
0&=&
(-1)^{m-1}\int e^{\frac i\hbar S_0}\prod_{i=1}^nS_{\alpha_i}
\frac i\hbar\left(
S_0,\mathcal{O}_x(f_0,\dots,f_m)
\right)
\\
&&+\int e^{\frac i\hbar S_0}
\sum_{1\leq j<k\leq n}\epsilon_{jk}
(S_{\alpha_j},S_{\alpha_k})\prod_{i\neq j,k}S_{\alpha_i}
\mathcal{O}_x(f_0,\dots,f_m).
\end{eqnarray*}
The sign is $\epsilon_{jk}=(-1)^{(g_1+\cdots+g_j)g_j+(g_1+\cdots
+g_{j-1}+g_{j+1}+\cdots+g_{k-1})g_k}$, where $g_j=p_j-1$ is the
ghost number of $S_{\alpha_j}$.
Now the BV bracket $(S_{\alpha_j},S_{\alpha_k})$ is again of the
form $S_\alpha$:
\[
(S_{\alpha_j},S_{\alpha_k})=-S_{[\alpha_j,\alpha_k]}.
\]
The Schouten--Nijenhuis bracket $[\ ,\ ]$ 
is a graded super Lie algebra structure on $T_{\mathrm{poly}}M$.
 On
vector fields it is defined to be the usual Lie bracket, and it is extended
to poly\-vector fields by the Leibnitz rule
 $[\alpha_1,\alpha_2\wedge\alpha_3]=[\alpha_1,\alpha_2]\wedge\alpha_3
+(-1)^{p_1(p_2-1)}\alpha_2\wedge[\alpha_1,\alpha_3]$. The Jacobi
identity for a bivector field $\alpha$ is $[\alpha,\alpha]=0$.

Moreover,
 $(S_0,f(\sX(t),\theta))=Df(\sX(t),\theta)$, 
which in components reads
\[
\left(S_0,f(X(t))\right)=0,\qquad
\left(S_0,\partial_if(X(t))\eta^{+i}(t)\right)=-d f(X(t)).
\]
Using this identity and the Leibnitz rule \Ref{e-leibnitz} we see that
the integral over
$B_m$ 
reduces to an integral over the boundary (of a
suitable compactification). The compactification may be understood
by thinking of $B_m$ as the moduli space of discs with $m+1$ marked
points on the boundary modulo the action of $SU(1,1)$.
The boundary of $\bar B_m$ consists of discs
degenerated into pairs of discs with a point in common. Its
various connected components are obtained by distributing
the points on the two discs in all possible ways compatible
with the cyclic ordering, see Fig.~2. 
\begin{figure}[h]
{\scalebox{.5}{\includegraphics{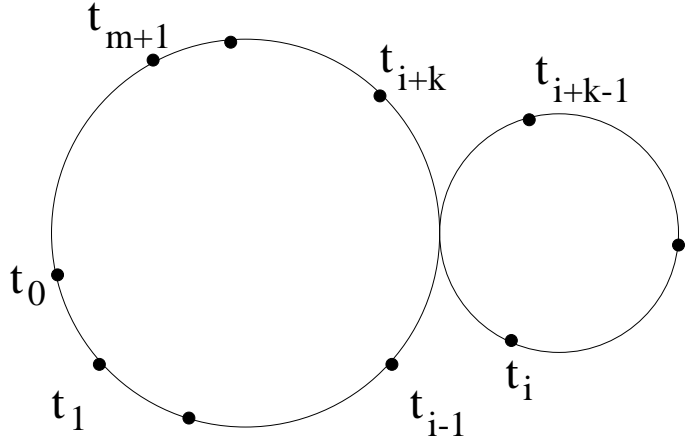}}}
\caption{A component of the boundary of $B_m$.  The point
$\infty$ is  $t_{m+1}$}
\end{figure}
By the usual factorization
arguments of quantum field theory, we obtain the following formulae.
Let $S_{\ell,n-\ell}$ be the subset of the group $S_n$ of permutations of
$n$ letters consisting of permutations such that $\sigma(1)<\cdots
<\sigma(\ell)$
and $\sigma(\ell+1)<\cdots<\sigma(n)$. 
For $\sigma\in S_{\ell,n-\ell}$ let us introduce the sign
\[
\epsilon(\sigma)=(-1)^{\sum_{r=1}^\ell g_{\sigma(r)}
(\sum_{s=1}^{\sigma(r)-1}g_s-\sum_{s=1}^{r-1}g_{\sigma(s)})}.
\]
It is the sign one gets if one
puts a product of $n$ elements of degree
$g_1$,\dots, $g_n$ of a graded commutative algebra in the order
given by $\sigma$. Then we have
\begin{eqnarray*}
\lefteqn{
\sum_{\ell=0}^n\sum_{k=1}^{m-1}
 \sum_{i=0}^{m-k}
\sum_{\sigma\in S_{\ell,n-\ell}}\epsilon(\sigma)(-1)^{k(i+1)}(-1)^m
U_{\ell}(\alpha_{\sigma(1)},\dots,\alpha_{\sigma(\ell)})
(f_0\otimes\cdots\otimes f_{i-1}}\\
&&
\otimes\, U_{n-\ell}(\alpha_{\sigma(\ell+1)},
\dots,\alpha_{\sigma(n)})(f_{i}\otimes\cdots\otimes f_{i+k})
\otimes f_{i+k+1}\otimes\cdots\otimes f_{m})
\\
&=&
\sum_{i<j}\epsilon_{ij} U_{n-1}([\alpha_i,\alpha_j],\alpha_1,\dots,
\widehat\alpha_i,\dots,\widehat\alpha_j,\dots,\alpha_n)
(f_0\otimes\cdots\otimes f_m).
\end{eqnarray*}
The sign $-(-1)^{k(i+1)}$ comes from the orientation of the
faces of the boundary of $\bar B_m$.
A sequence of maps $U_n$ with this property and the symmetry
\Ref{e-sym} is called an
$L_\infty$ morphism \cite{SchlSt, LS, K}. A consequence of it in this case
is  the
 formality conjecture: $U_1$, which is (up to sign)
the obvious map sending a poly\-vector
field to itself viewed as poly\-differential operator, induces
an isomorphism of graded Lie algebras from the graded Lie algebra
of poly\-vector fields to the cohomology of the poly\-differential
operators, viewed as a complex of
the Hochschild cochains of $C^\infty(M)$, see \cite{K}.

A special case of this identity is the associativity of
the star product: in this case $\alpha$ is a bivector field
(a poly\-vector field of degree $1$) obeying the Jacobi identity
$[\alpha,\alpha]=0$. Then we have $U(\alpha)=\sum_{n=0}^\infty
\hbar^nU_n(\alpha,\dots,\alpha)$ and by the ghost number
condition \Ref{e-gh}, every $U_n$ is a bidifferential operator.
The $L_\infty$ identity reduces to the associativity of $U(\alpha)$:
\[
\sum_{k=0}^nU_{k}(U_{n-k}(f_0\otimes f_1)\otimes f_2)
-
\sum_{k=0}^nU_{k}(f_0\otimes U_{n-k}(f_1\otimes f_2))=0,
\]
where we have suppressed the dependence on $\alpha$ in the
notation.
\subsection{Tadpoles}
The perturbative expansion of $U(\alpha)$ contains 
tadpoles which can be removed as in \ref{ss-ren} either
by choosing a constant angle $\vartheta$, or by replacing
$S_\alpha$ by 
\[
S_\alpha'=S_\alpha-\frac{i\hbar}{2\pi}
\frac1{p!}\int_D\int d^2\theta\,
\tilde\kappa
\partial_k\alpha^{k,j_1,\dots,j_{p}}
\seta_{j_1}\cdots\seta_{j_{p}}.
\]
The arguments of the previous subsections remain valid
with the addition of these counterterms, since
$\triangle S_\alpha'=(S_0,S_\alpha')=0$ and
$(S_\alpha',S_\beta')=-S_{[\alpha,\beta]}'$.


\begin{thebibliography}{BFFLS}
\bibitem[AKSZ]{AKSZ} M. Alexandrov, M. Kontsevich, A. Schwarz and O.
Zaboronsky, 
{\it The geometry of the master equation and topological quantum
field theory}, Internat. J. Modern Phys. A 12 (1997), no.\ 7, 1405--1429
\bibitem[BV]{BV} I. Batalin and G. Vilkovisky,
{\it Gauge algebra and quantization},
Phys.\ Lett.\ 102 B (1981), 27;
{\it
Quantization of gauge theories with linearly dependent
generators}, Phys.\ Rev.\ D29 (1983), 2567
\bibitem[BT]{BT} M. Blau, G. Thompson,
{\it Topological gauge theories of antisymmetric 
tensor fields}, Ann. Physics 205 (1991), no. 1,
130--172
\bibitem[BFFLS]{BFFLS}
F. Bayen, M. Flato, C. Fronsdal, A. Lichnerowicz and D. Sternheimer,
{\it Deformation theory and quantization I, II}, 
Ann. Phys.  111 (1978), 61--110, 111--151
\bibitem[CFT]{CFT}
A. S. Cattaneo, G. Felder and L. Tomassini,
in preparation
\bibitem[DDV]{DVV} R. Dijkgraaf, H. Verlinde and E. Verlinde,
{\it Topological strings in $d<1$},
Nucl. Phys. B 352 (1991), 59--86
\bibitem[F]{F} R. P. Feynman,
{\it Space--time approach to non-relativistic quantum mechanics},
Rev.\ Modern Physics 20, (1948), 367--387
\bibitem[I]{I}
N. Ikeda,
{\it
Two-dimensional gravity and nonlinear gauge theory},
Ann. Phys. 235, (1994) 435--464
\bibitem[K]{K} M. Kontsevich,
{\it Deformation quantization of Poisson manifolds},
q-alg/9709040
\bibitem[LS]{LS} T. Lada and J. Stasheff,
{\it Introduction to sh Lie algebras for physicists},
Internat. J. Theoret. Phys. 32 (1993), no. 7, 1087--1103
\bibitem[SchStr]{SchStr}
P. Schaller and T. Strobl,
{\it  Poisson structure induced (topological) field theories}
Modern Phys. Lett. A 9 (1994), no. 33,
3129--3136
\bibitem[SchlSt]{SchlSt} M. Schlessinger and J. Stasheff,
{\it The Lie algebra structure of tangent cohomology and
deformation theory},
J. Pure Appl. Alg. 38 (1985), 313--322
\bibitem[S1]{S1} A. S. Schwarz,
{\it Geometry of Batalin--Vilkovisky quantization},
Commun. Math. Phys. 155 (1993), 249--260
\bibitem[S2]{S2} A. S. Schwarz,
{\it The partition function of a degenerate quadratic functional and
the Ray--Singer invariants}, Lett. Math. Phys. 2 (1978), 247--252
\bibitem[W1]{W1} E. Witten, 
{\it A note on the antibracket formalism},
Modern Phys. Lett. A 5 (1990), no.\ 7, 487--494
\bibitem[W2]{W2} E. Witten,
{\it On the structure of the topological phase of two-dimensional gravity},
Nucl. Phys. B 340 (1990) 281--332
\end{thebibliography}
\end{document}